\documentclass[12pt]{amsart}

\usepackage{latexsym}
\usepackage{amssymb}
\usepackage{amsmath}
\usepackage{amscd}
\usepackage{mathrsfs}

\newtheorem{theo}{Theorem}[section]
\newtheorem{lemma}[theo]{Lemma}
\newtheorem{prop}[theo]{Proposition}
\newtheorem{cor}[theo]{Corollary}

\newtheorem{quest}{Question}[section]

\theoremstyle{definition}

\newtheorem{rem}{Remark}[section]

\DeclareMathOperator{\Der}{{\rm D}}
\DeclareMathOperator{\ST}{{{\rm T}^1}}
\DeclareMathOperator{\Vis}{{\rm Vis}}

\DeclareMathOperator{\Span}{{\rm span}}
\DeclareMathOperator{\GL}{{\rm GL}}
\DeclareMathOperator{\SL}{{\rm SL}}
\DeclareMathOperator{\SO}{{\rm SO}}
\DeclareMathOperator{\End}{{\rm End}}

\DeclareMathOperator{\nor}{{\rm N}_G} 
\DeclareMathOperator{\norone}{{\rm N}^1_G}
\DeclareMathOperator{\cen}{{\rm Z}_G} 
\DeclareMathOperator{\cenKA}{K_0}           
\DeclareMathOperator{\supp}{{\rm supp}} 
\DeclareMathOperator{\Stab}{{\rm Stab}} 
\DeclareMathOperator{\Ad}{{\rm Ad}} 
\DeclareMathOperator{\diag}{\rm{diag}} 
\DeclareMathOperator{\Cc}{{\rm C}_{{\rm c}}} 
\DeclareMathOperator{\Ctwo}{{\rm C}^2} 
\DeclareMathOperator{\dd}{{\rm d}}

\renewcommand{\phi}{\varphi}

\newcommand{\del}{\partial}

\newcommand{\field}[1]{\mathbb{#1}} 
 
\newcommand{\R}{\field{R}}
\newcommand{\N}{\field{N}}

\newcommand{\Z}{\field{Z}}
\newcommand{\HH}{\field{H}}
\newcommand{\T}{\field{T}}
\newcommand{\Sn}{\field{S}}

\providecommand{\abs}[1]{\lvert#1\rvert}
\providecommand{\Abs}[1]{\left\lvert #1 \right\rvert}
\providecommand{\norm}[1]{\lVert#1\rVert}

\renewcommand{\setminus}{\smallsetminus}

\newcommand{\inv}{^{-1}}

\newcommand{\cl}[1]{\overline{#1}}

\newcommand{\tcup}[1]{\textstyle{\bigcup_{#1}}\,}

\newcommand{\psmat}[1]{\bigl(\begin{smallmatrix} #1
\end{smallmatrix}\bigr)}

\newcommand{\la}[1]{\mathfrak{\lowercase{#1}}}

\newcommand{\cA}{\mathcal{A}}

\newcommand{\cC}{\mathcal{C}}

\newcommand{\cD}{\mathcal{D}}

\newcommand{\cE}{\mathcal{E}}

\newcommand{\cF}{\mathcal{F}}
\newcommand{\sF}{\mathscr{F}}

\newcommand{\cH}{\mathcal{H}}
\newcommand{\sH}{\mathscr{H}}

\newcommand{\cJ}{\mathcal{J}}

\newcommand{\cK}{\mathcal{K}}

\newcommand{\cO}{\mathcal{O}}

\newcommand{\cS}{\mathcal{S}}

\newcommand{\gmg}{{G/\Gamma}}
\newcommand{\nbd}{be a neighbourhood of $e$ in $G$ such that}

\begin{document}

\title[Limiting distributions of curves under geodesic flow]{Limiting
  distributions of curves under geodesic flow on hyperbolic manifolds}

\author{Nimish A.  Shah} \address{Tata Institute of Fundamental
  Research, Mumbai 400005, INDIA} 
\date{March 21, 2008}
\email{nimish@math.tifr.res.in}

\begin{abstract}
  We consider the evolution of a compact segment of an analytic curve
  on the unit tangent bundle of a finite volume hyperbolic
  $n$-manifold under the geodesic flow. Suppose that the curve is not
  contained in a stable leaf of the flow. It is shown that under the
  geodesic flow, the normalized parameter measure on the curve gets
  asymptotically equidistributed with respect to the normalized
  natural Riemannian measure on the unit tangent bundle of a closed
  totally geodesically immersed submanifold.

  Moreover, if this immersed submanifold is a proper subset, then a
  lift of the curve to the universal covering space $\ST(\HH^n)$ is
  mapped into a proper subsphere of the ideal boundary sphere
  $\del\HH^n$ under the visual map.  This proper subsphere can be
  realized as the ideal boundary of an isometrically embedded
  hyperbolic subspace in $\HH^n$ covering the closed immersed
  submanifold.

  In particular, if the visual map does not send a lift of the curve
  into a proper subsphere of $\del\HH^n$, then under the geodesic flow
  the curve gets asymptotically equidistributed on the unit tangent
  bundle of the manifold with respect to the normalized natural
  Riemannian measure.

  The proof uses dynamical properties of unipotent flows on finite
  volume homogeneous spaces of $\SO(n,1)$. 
\end{abstract}

\maketitle

\section{Introduction}

It is instructive to note the following dynamical property: Let
$\psi:I=[0,1]\to \R^n$ be a $\Ctwo$-curve such that for any proper
rational hyperplane, say $\cH$, in $\R^n$, the set $\{s\in
I:\psi(s)\in\cH\}$ has null measure. Let $\T^n=\R^n/\Z^n$, and let
$\pi:\R^n\to\T^n$ denote the quotient map.  Then for any continuous
function $f$ on $\T^n$,
\begin{equation}
  \label{eq:Tn} 
  \lim_{\alpha\to\infty} \int_0^1 f(\pi(\alpha\psi(s)))\,{\dd} s = 
  \int_{\T^n} f(x)\,{\dd}x,
\end{equation}
where ${\dd}x$ denotes the normalized Haar integral on $\T^n$.  Using
Fourier transforms, we can verify \eqref{eq:Tn} for the characters
\[
f_m(x):=\exp(2\pi(m\cdot x)),\quad \textrm{$\forall x\in \T^n$, where
  $m\in\Z^n$}.
\]

The above observation was used in \cite{Lin+Klaus:Invariant} for
$\psi(t)=(\cos(2\pi t),\sin(2\pi t))$; the unit circle in $\R^2$. Later
we learnt that a general result in this direction was obtained
earlier by B. Randol~\cite{Randol:dilate} in response to a question
raised by D. Sullivan.

Now we ask a similar question for the hyperbolic spaces. Consider the
unit ball model $B^n$ for the hyperbolic $n$-space $\HH^n$ of constant
curvature $(-1)$.  Let $\Gamma\subset \SO(n,1)$ be a discrete subgroup
such that $M:=\HH^n/\Gamma$ is a hyperbolic manifold of finite
Riemannian volume.  Let $\pi:\HH^n\to M$ be the quotient map.  As a
special case of a more general result proved in
\cite{D+R+S:symm,Esk+McM:mixing}, we have that if we project the
invariant probability measure on the sphere $\alpha \Sn^{n-1}\subset
B^n$, for $0<\alpha<1$, under $\pi$ to $M$, then asymptotically as
$\alpha\to 1^-$, the measure gets equidistributed with respect to the
normalized measure associated to the Riemannian volume form on
$M$. The case of $n=3$ was proved earlier in \cite{Randol:dilate}. 

In this article, we will address the following much more refined
problem: Instead of the invariant measure on the sphere, we take a
smooth measure on a one-dimensional curve on $\Sn^{n-1}$ and describe
the limiting distribution of the projection of its expands on
$\alpha\Sn^{n-1}$ as $\alpha\to 1^-$.

\begin{theo}
  \label{thm:Hn-basic}
  Let $\bar\psi:I=[0,1]\to \Sn^{n-1}$ be an analytic map.  If
  $\bar\psi(I)$ is not contained in a proper subsphere in $\Sn^{n-1}$,
  then for any $f\in\Cc(M)$,
  \begin{equation}
    \label{eq:Hn}
    \lim_{\alpha\to 1^-} \int_If(\pi(\alpha\bar\psi(s))\,{\dd}s=\int_Mf(x)\,{\dd}x,
  \end{equation}
  where ${\dd}x$ denotes the normalized integral associated to the
  Riemannian volume form on $M$.
\end{theo}

By a {\em proper subsphere\/} of $\Sn^{n-1}\subset\R^n$ we mean the
intersection of $\Sn^{n-1}$ with a proper affine subspace of $\R^n$.

Now we describe a generalization of the above phenomenon in a suitable
geometric framework.  Let $\del\HH^n$ denote the ideal boundary of
$\HH^n$. Let $\ST(\HH^n)$ denote the unit tangent bundle on $\HH^n$.
We identify $\del\HH^n$ with $\Sn^{n-1}$.  Let
\[
\Vis:\ST(\HH^n)\to \del\HH^n\cong \Sn^{n-1},
\]
denote the visual map sending a tangent to the equivalence class of
the directed geodesics tangent to it. Thus any fiber of the visual
map is a (weakly) stable leaf of the geodesic flow.  Now let
$M$ be any $n$-dimensional hyperbolic manifold (with constant
curvature $(-1)$ and) with finite Riemannian volume, let $\ST(M)$
denote the unit tangent bundle on $M$, and let $\{g_t\}$ denote the
geodesic flow on $\ST(M)$.  Let $\pi:\HH^n\to M$ be a universal
covering map, and let $\Der\pi:\ST(\HH^n)\to \ST(M)$ denote its
derivative.

\begin{theo}
  \label{thm:Hn}
  Let $\psi:I=[a,b]\to \ST(M)$ be an analytic curve such that
  $\Vis(\tilde\psi(I))$ is not a singleton set, where
  $\tilde\psi:I\to\ST(\HH^n)$ denotes a lift of $\psi$ to the covering
  space; that is, $\Der\pi\circ\tilde\psi=\psi$.  Then there exists a
  totally geodesic immersion $\Phi:M_1\to M$ of a hyperbolic manifold
  $M_1$ with finite volume such that the following holds: $\forall
  f\in\Cc(\ST(M))$,
  \begin{equation}
    \label{eq:M1-measure}
    \lim_{t\to\infty} \frac{1}{\abs{I}}\int_I f(g_{t}\psi(s))\,{\dd}s 
    = \int_{\ST(M_1)}f((\Der\Phi)(v))\,{\dd}v,
  \end{equation}
  where $\abs{\cdot}$ denotes the Lebesgue measure, and ${\dd}v$
  denotes the normalized integral on $\ST(M_1)$ associated to the
  Riemannian volume form on $M_1$.

  Moreover if $\pi':\HH^m\to M_1$ denotes a locally isometric covering
  map, then there exists an isometric embedding
  $\tilde\Phi:\HH^m\hookrightarrow\HH^n$ such that
  \[
  \pi\circ\tilde\Phi=\Phi\circ\pi' \text{ and }
  \Vis(\tilde\psi(I))\subset\del(\tilde\Phi(\HH^m)).
  \]
\end{theo}

In order to describe the relation between $\Vis(\tilde\psi(I))$ and
the totally geodesic immersion $\Phi$, we will recall the following:

\begin{theo}[\cite{R:uniform},\cite{Shah:tot-geod}]
  \label{thm:tot-geod}
  Let $M$ be a hyperbolic manifold with finite Riemannian volume.  For
  $k\geq 2$, let $\Psi:\HH^k\to M$ be a totally geodesic immersion.
  Then there exists a totally geodesic immersion $\Phi:M_1\to M$ of a
  hyperbolic manifold $M_1$ with finite Riemannian volume such that
  \[
  \cl{\Psi(\HH^k)}=\Phi(M_1) \quad \textrm{and} \quad
  \cl{\Der\Psi(\ST(\HH^k))}=\Der\Phi(\ST(M_1)).
  \]
  \qed
\end{theo}

This result can be obtained as a direct consequence of the orbit
closure theorem for unipotent flows (Raghunathan's conjecture) proved by
Ratner~\cite{R:uniform}; more specifically, the fact that the
closure of any $\SO(k,1)$-orbit in $\SO(n,1)/\Gamma$ is a closed orbit
of a subgroup of the form $Z\cdot\SO(m,1)$, where $Z$ is a compact
subgroup of the centralizer of $\SO(m,1)$ in $\SO(n,1)$.

\begin{rem}
  \label{rem:M1}
  Let the notation be as in Theorem~\ref{thm:Hn}.  Let $\Sn^{k-1}$ be
  the smallest dimensional subsphere of $\del\HH^n\cong \Sn^{n-1}$
  such that $\Vis(\tilde\psi(I))\subset \Sn^{k-1}$.  Since
  $\Vis(\tilde\psi(I))$ is not a singleton set, we have $2\leq k\leq
  n$.  Therefore there exists an isometric embedding
  $\HH^k\hookrightarrow \HH^n$ such that $\del\HH^k=\Sn^{k-1}$.  If
  $\{\tilde g_t\}$ denotes the geodesic flow and $\tilde
  d(\cdot,\cdot)$ denotes the distance function on $\ST(\HH^n)$, then
  \begin{equation}
    \label{eq:91}
    \lim_{t\to\infty} \sup_{s\in I} \tilde d(\tilde g_t\tilde\psi(s),\ST(\HH^k))=0.
  \end{equation}

  Since $\pi:\HH^k\to M$ is a totally geodesic immersion, by
  Theorem~\ref{thm:tot-geod} there exists a totally geodesic immersion
  $\Phi:M_1\to M$ of a hyperbolic manifold of finite Riemannian volume
  such that
  \begin{equation}
    \label{eq:M1}
    \Phi(M_1)=\cl{\pi(\HH^k)}.
  \end{equation}
  This describes the map $\Phi$ as involved in the statement of
  Theorem~\ref{thm:Hn}.  Also, by \eqref{eq:91} and \eqref{eq:M1}, if
  $d(\cdot,\cdot)$ denotes the distance function on $M$ then
  \begin{equation}
    \label{eq:92} 
    \lim_{t\to\infty} \sup_{s\in I}d(g_t\psi(s),\Der\Phi(\ST(M_1)))=0.
  \end{equation}
\end{rem}

We have the following consequences.

\begin{theo}
  \label{thm:S-2}
  Let $M$ be a hyperbolic Riemannian manifold with finite volume.  Let
  $\psi:I\to M$ be an analytic map such that $\Vis(\tilde\psi(I))$ is
  not contained in a proper subsphere in $\del\HH^n$, where
  $\tilde\psi:I\to \ST(\HH^n)$ is a lift of $\psi$ such that
  $\Der\pi\circ\tilde\psi=\psi$.  Then given any $f\in\Cc(\ST(M))$,
  \[
  \lim_{t\to\infty} \frac{1}{\abs{I}}\int_I
  f(a_t\psi(s))\,{\dd}s=\int_{\ST(M)} f\,{\dd}v,
  \]
  where ${\dd}v$ is the normalized integral on $\ST(M)$ associated to
  the Riemannian volume form on $M$.
\end{theo}

\begin{cor}
  \label{cor:S}
  Let $M$ be a hyperbolic manifold with finite volume.  Let $x\in M$
  and $\psi:I=[a,b]\to \ST_x(M)$ be an analytic map such that
  $\psi(I)$ is not contained in any proper subsphere in $\ST_x(M)$.
  Then
  \[
  \lim_{t\to\infty} \frac{1}{\abs{I}}\int_I
  f(g_t\psi(s))\,{\dd}t=\int_{\ST(M)}f(v)\,{\dd}v,\quad\forall
  f\in\Cc(\ST(M)),
  \]
  where ${\dd}v$ is the normalized Riemannian volume integral on
  $\ST(M)$.
\end{cor}

It may be interesting to compare the above result with
\cite{Zeghib:invariant-I} where any rectifiable invariant set for the
geodesic flow is shown to be a conull subset of the unit tangent
bundle of a closed finite volume totally geodesic submanifold.

\subsection{Reformulation in terms of flows on homogeneous spaces}
 
Let $G=\SO(n,1)$, and $P^-$ be a minimal parabolic subgroup of $G$,
and $K\cong\SO(n)$ be a maximal compact subgroup of $G$.  Then
$M:=P^-\cap K\cong \SO(n-1)$.  Since $G=P^-K$,
\begin{equation}
  \label{eq:SOn}
  P^-\backslash G\cong M\backslash K\cong\SO(n-1)\backslash SO(n)\cong
  \Sn^{n-1}.
\end{equation}
We let $p:G\to \Sn^{n-1}$ denote the quotient map corresponding to
\eqref{eq:SOn}.  Let $A$ be a maximal connected $\R$-diagonalizable
subgroup of $G$ contained in $\cen(M)\cap P^-$.  Since $G$ is of
$\R$-rank $1$, $A$ is a one-parameter group, and the centralizer of
$A$ in $G$ is $\cen(A):=MA$.  Let $N^-$ denote the unipotent radical
of $P^-$.  Define
\begin{eqnarray}
  A^+ &=&\{a\in A: a^k g a^{-k}\to e \textrm{ as } k\to\infty\textrm{ for
    any }g\in N^-\}, \textrm{ and} \\
\label{eq:N}
  N & =&\{g\in G: a^{-k} g a^k\to e \textrm{ as } k\to\infty \textrm{ for
    any }a\in A^+\}.
\end{eqnarray}
Let $\la{n}$ denote the Lie algebra on $N$.  Then $\la{n}$ is abelian,
and we identify it with $\R^{n-1}$.  Let $u:\R^{n-1}\to N$ be the map
$u(v)=\exp(v)$ for any $v\in\R^{n-1}\cong \la{n}$. We observe that the
map
\begin{equation}
  \label{eq:stereo}
  S:\R^{n-1}\to \Sn^{n-1} \textrm{ defined by } S(v)=p(u(v)), \quad
  \forall v\in \R^{n-1}, 
\end{equation}
is the inverse of stereographic projection.

Let $\alpha:A\to\R^\ast$ be the character such that
$au(v)a\inv=u(\alpha(a)v)$ for all $v\in\R^{n-1}$.  Then $A^+=\{a\in
A:\alpha(a)>1\}$.

Let $\Gamma$ be a lattice in $G$ and $\mu_G$ be the $G$-invariant
probability measure on $\gmg$.

\begin{theo}
  \label{thm:main}
  Let $\theta:I=[a,b]\to G$ be an analytic map such that
  $p(\theta(I))$ is not contained in a subsphere of $\Sn^{n-1}$.  Then
  given any $f\in \Cc(\gmg)$, any compact set $\cK\subset G/\Gamma$
  and any $\epsilon>0$, there exists $R>0$ such that for any $a\in
  A^+$ with $\alpha(a)>R$,
  \begin{equation}
    \label{eq:mainthm}
    \Abs{\frac{1}{\abs{I}}\int_I f(a\theta(t)x))\,{\dd}t 
      - \int_{\gmg}f\,{\dd}\mu_G}<\epsilon,\quad \forall x\in\cK.   
  \end{equation}
\end{theo}

First we shall consider the following crucial case of the above
theorem.
 
\begin{theo}
  \label{thm:main-unip}
  Let $\phi:I=[a,b]\to\R^{n-1}$ be an analytic curve such that
  $\phi(I)$ is not contained in any sphere or an affine hyperplane.
  Let $x_i\stackrel{i\to\infty}{\longrightarrow} x$ be a convergent
  sequence in $\gmg$, and let $\{a_i\}_{i\in\N}$ be a sequence in
  $A^+$ such that
  $\alpha(a_i)\stackrel{i\to\infty}{\longrightarrow}\infty$.  Then
  \begin{equation}
    \label{eq:main-unip}
    \lim_{i\to\infty} \frac{1}{\abs{I}}\int_I
    f(a_iu(\phi(t)x_i)\,{\dd}t=\int_\gmg f\,{\dd}\mu_G,\quad\forall
    f\in\Cc(\gmg).
  \end{equation}
\end{theo}

We will deduce the above result from the following general statement,
which is the main result of this paper.

\begin{theo}
  \label{thm:psi}
  Let $\phi:I\to \R^{n-1}$ be a nonconstant analytic map, and let
  $x\in G/\Gamma$. Then there exist a closed subgroup $H$ of $G$, an
  analytic map $\zeta:I\to M(=\cen(A)\cap K)$ and $h_1\in G$ such that
  $\pi(H)$ is closed and admits a finite $H$-invariant measure, say
  $\mu_H$, and the following holds: For any sequence
  $\{a_i\}_{i\in\N}\subset A^+$, if
  $\alpha(a_i)\stackrel{i\to\infty}{\longrightarrow}\infty$ then
  \begin{equation}
    \label{eq:lambda-muH} 
    \lim_{i\to\infty} \int_If(a_iu(\phi(t))x)\,{\dd}t 
    = \int_{t\in I}\left(\int_{y\in\gmg}f(\zeta(t)h_1y)\,{\dd}\mu_H\right){\dd}t.
  \end{equation}
  Moreover $A\subset h_1Hh_1\inv$, $N\cap h_1Hh_1\inv\neq\{e\}$, and
  there exists $g\in G$ such that $x=\pi(g)$ and
  \begin{equation}
    \label{eq:UHg} 
    u(\phi(t))g\in N^-\zeta(t)h_1H,\quad \forall t\in I.
  \end{equation}
\end{theo}

\begin{rem} 
  Suppose we are given a convergent sequence $x_i\to x$ in
  $G/\Gamma$. We consider \eqref{eq:lambda-muH} for $x_i$ in place of
  $x$ in the statement of Theorem~\ref{thm:psi}. Then the limiting
  distribution depends on the choice of the sequence $\{a_i\}$. We can
  still conclude that the analogue of \eqref{eq:lambda-muH} holds
  after passing to a subsequence.
\end{rem}

\subsubsection*{Acknowledgment} {\small The precise geometric
  consequences of the equidistribution results on homogeneous spaces
  were formulated during author's visits to Elon Lindenstrauss at
  Princeton University and Hee Oh at Brown University. The author
  would like to thank them for their hospitality and support. Special
  thanks are due to the referee for the simplified proof of
  Theorem~\ref{thm:W-invariant} given here, and many other useful
  suggestions.}

\section{Non-divergence of translated measures}
\label{sec:nondiv}
Let $\phi:I\to \R^{n-1}$ be a nonconstant analytic map.  Let
$\{a_i\}\subset A^+$ be a sequence such that
$\alpha(a_i)\stackrel{i\to\infty}{\longrightarrow}\infty$.  Let
$x_i\stackrel{i\to\infty}{\longrightarrow} x$ be a convergent sequence
in $\gmg$.  For each $i\in\N$, let $\mu_i$ be the measure on $\gmg$
defined by
\begin{equation}
  \label{eq:mui} 
  \int_{\gmg} f\,{\dd}\mu_i:=\frac{1}{\abs{I}}\int_I
  f(a_iu(\phi(t))x_i)\,{\dd}t\quad\forall f\in\Cc(\gmg).
\end{equation}

This section is devoted to the proof of the following:

\begin{theo}
  \label{thm:return}
  Given $\epsilon>0$ there exists compact set $\cK\subset G/\Gamma$
  such that $\mu_i(\cK)\geq 1-\epsilon$ for all $i\in\N$.
\end{theo}

We begin with some notation.  Let $d=\dim N$, $\la{g}$ denote the Lie
algebra of $G$, and $V=\wedge^d\la{g}$.  Let
$\boldsymbol{p}\in\wedge^d\la{N}\setminus\{0\}$.  Consider the
$\wedge^d \Ad$-action of $G$ on $V$.

\subsection{$(C,\alpha)$-good family}
\label{subsec:CAlpha-1}
Let $\sF$ be the $\R$-$\Span$ of the coordinate functions of the map
$\Upsilon:I\to\End(V)$ given by $\Upsilon(t)=\wedge^d\Ad(u(\phi(t)))$
for all $t\in I$.

Fix $t_0\in I$ and let $\cE$ be the smallest subspace of $\End(V)$
such that $\Upsilon(I)\subset \cE+\Upsilon(t_0)$. Then
$\Upsilon(I)\subset \cE+\Upsilon(t)$ for all $t\in I$. For any $t\in
I$, we have $\cE_t:=\Span\{\Upsilon^{(k)}(t):k\geq 1\}\subset \cE$,
where $\Upsilon^{(k)}(t)$ denotes the $k$-th derivative at $t$.  Since
$\Upsilon$ is an analytic function, we have $\Upsilon(I)\subset
\Upsilon(s)+\cE_s$. Therefore $\cE\subset \cE_s$. Hence $\cE_s=\cE$
for all $s\in I$.
 
Therefore by \cite[Proposition~3.4]{Klein+Mar:Annals98}, applied to
the function $t\mapsto \Upsilon(t)-\Upsilon(t_0)$ from $I$ to $\cE$,
there exist constants $C>0$ and $\alpha>0$ such that the family $\sF$
consists of {\em $(C,\alpha)$-good functions\/}; that is, for any
subinterval $J\subset I$, $\xi\in\cF$ and $r>0$,
\begin{equation} 
\label{eq:CAlpha} 
\abs{\{t\in J: \abs{\xi(t)}<r\}} \leq 
C\left(\frac{r}{\sup_{t\in J}|\xi(t)|}\right)^\alpha\abs{J}.
\end{equation}
It may be noted that, since $I$ is compact, by the result quoted
above, a priori \eqref{eq:CAlpha} holds only for subintervals $J$ with
$\abs{J}$ smaller than a fixed constant depending on $\Upsilon$ and
$I$. Then by a straightforward argument using a finite well-overlapping
covering of $I$ by short intervals of fixed length, and applying the
above inequality successively, we can choose a much larger $C$ such
that \eqref{eq:CAlpha} for all subintervals $J\subset I$.

Now we fix a norm $\norm{\cdot}$ on $V$.  Then given any $\epsilon>0$
  and $r>0$, there exists $R>0$ such that for any $h_1,h_2\in G$ and
  an interval $J\subset I$, one of the following holds:
\begin{enumerate}
\item[I)] $\sup_{t\in J}\norm{h_1u(\phi(t))h_2\boldsymbol{p}}<R$.
\item[II)] 
  \[
\Abs{\{t\in J: \norm{h_1u(\phi(t))h_2 \boldsymbol{p}}\leq r\}} \leq \epsilon 
   \Abs{\{t\in J:\norm{h_1u(\phi(t))h_2 \boldsymbol{p}}\leq R\}}
\].
\end{enumerate}

\begin{prop}[\cite{Dani:rk=1}]
  \label{prop:return} 
  There exists a finite set $\Sigma\subset G$ such that $\Gamma \Sigma
  \boldsymbol{p}$ is a discrete subset of $V$, and the following
  holds: Given $\epsilon>0$ and $R>0$, there exists a compact set
  $\cK\subset\gmg$ and such that for any $h_1,h_2\in G$, and a
  subinterval $J\subset I$, one of the following holds:
  \begin{enumerate}
  \item[I)] 
    There exists $\gamma\in\Gamma$ and $\sigma\in\Sigma$ such that
    \begin{equation*}
      \sup_{t\in J} \norm{h_1u(\phi(t))h_2\sigma\gamma \boldsymbol{p}}<R.   
    \end{equation*}
  
  \item[II)]
      $\abs{t\in J:\pi(h_1u(\phi(t))h_2)\in\cK\}}\geq (1-\epsilon)\abs{J}$.
  \end{enumerate}
\end{prop}

In the above proposition, $(\sigma N \sigma\inv)\cap \Gamma$ is a
cocompact lattice in $\sigma N\sigma\inv$ for each $\sigma\in\Sigma$.

Now we will make an observation which will allow us to prove that the
possibility~(I) in the conclusion of the above proposition will not
hold in the situation of our interest.
 
\subsection{Basic lemma}

Consider a linear representation of $\SL(2,R)$ on a finite
dimensional vector space $V$. Let $a=\psmat{\alpha & \\ & \alpha\inv}$
for some $\alpha>1$, and define
\begin{align}
\begin{split}
  V^+&=\{v\in V: a^{-k}v\to\infty\ \textrm{as $k\to\infty$}\} \\
  V^0& =\{v\in V: av=v\}\\
  V^-& =\{v\in V: a^kv\to 0 \ \textrm{as $k\to\infty$}\}.
\end{split}
\end{align}
Then any $v\in V$ can be uniquely expressed as $v=v^++v^0+v^-$, where
$v^\pm\in V^\pm$ and $v^0\in V^0$.  We also write $V^{+0}=V^++V^0$,
and $V^{0-}=V^0+V^-$.  Let $q^+:V\to V^+$, $q^0:V\to V^0$,
$q^{+0}:V\to V^{+0}$, and $q^{0-}:V\to V^{0-}$ denote the projections
$q^+(v)=v^+$, $q^0(v)=v^0$, $q^{+0}(v)=v^{+0}:=v^++v^0$, and
$q^{0-}(v)=v^{0-}:=v^0+v^-$ for all $v\in V$.  We consider the
Euclidean norm on $V$ such that $V^+$, $V^0$ and $V^-$ are orthogonal.

\begin{lemma}
  \label{lem:SL2}
  Let $\boldsymbol{u}=\psmat{1 & t \\0 & 1}$ for some $t\neq 0$.  Then
  there exists a constant $\kappa=\kappa(t)>0$ such that
  \begin{equation} 
    \label{eq:abs1}
    \max\{\norm{v^{+}},\norm{(\boldsymbol{u}v)^{+0}}\} 
    \geq \kappa\norm{v},\quad \forall v\in V.
  \end{equation}
\end{lemma}

\begin{proof}
Since it is enough to prove the result for each of the
$\SL(2,\R)$-irreducible subspace of $V$. Therefore without loss of
generality we may assume that $\SL(2,\R)$ acts irreducibly on $V$.
 
  Let $m=\dim V -1$.  Then $m=2r-1$ or $m=2r$ for some
  $r\in\N$.  Consider the associated representation of the Lie algebra
  $\la{SL}_2(\R)$ on $V$.   Let $\boldsymbol{e}=\psmat{0&1\\0&0}$,
  $\boldsymbol{h}=\psmat{1& \\ &-1}$, and
  $\boldsymbol{f}=\psmat{0&0\\1&0}$ denote the standard
  $\la{SL}_2$-triple.  Then there exists a basis of $V$ consisting of
  elements $v_0,v_1,\dots,v_m$ such that
  \[
  \boldsymbol{h}v_k=(m-2k)v_k,\quad \textrm{and}\quad
  \boldsymbol{e}v_k=kv_{k-1}, \quad \forall 0\leq k\leq m,
  \]
  where $v_{-1}=0$.  Then
  \begin{equation}
    \label{eq:V0-pm}  
    V^{+0}=\Span\{v_0,\dots,v_{m-r}\} \quad \textrm{and} \quad
    V^{0-}=\Span\{v_r,\dots,v_m\}.
  \end{equation}
  Since $\boldsymbol{u}=\exp(t\boldsymbol{e})$, we have
  \[
  \boldsymbol{u}v_k=\sum_{l=0}^k \binom{k}{l} t^{k-l} v_l, \quad 0\leq
  k\leq m.
  \]
  Let $A$ denote the restriction of the map $\boldsymbol{u}$ from
  $V^+$ to $V^+$ with respect to the basis
  $\{v_0,\dots,v_{r-1}\}$.  Let $B$ denote the matrix of the map
  $q^{+0}\circ \boldsymbol{u}:V^{0-}\to V^{+0}$ with respect to the
  basis given by \eqref{eq:V0-pm}.

  Next we want to show that $B$ is invertible.   We write
  $b_{k,l}=t^{k-l}\binom{k}{l}$ for $r\leq k\leq m$ and $0\leq l \leq
  m-r$.   And for any $r\leq m_1\leq m$, we consider the
  $(m_1-r+1)\times (m_1-r+1)$-matrix
  \[
  B(m_1,r)=\bigl(b_{k,l}\bigr)_{\substack{r\leq k\leq m_1 \\ 0\leq l \leq m_1-r}}.
  \]
  Then $B=B(m,r)$.  In view of the binomial relations
  \[
  \binom{k+1}{l+1}-\binom{k}{l+1}=\binom{k}{l} \quad \textrm{and}
  \quad b_{k+1,l+1}-tb_{k,l+1}=b_{k,l},
  \]
  we apply the row operations $R_{k+1}-tR_k$, successively, in the
  order $k=(m_1-1),\dots,1$.  We obtain that
  \[
  \det B(m_1,r)=t^r\det B(m_1-1,r).
  \]
  Since $\det B(r,r)=t^r$, we get
  \[
  \det B=\det B(m,r)=t^{r(m-r+1)}.
  \]
  Since $t\neq 0$, $B$ is invertible.

  Now
  \begin{equation} 
    \label{eq:abs2}
    \norm{(\boldsymbol{u}v)^{+0}}=\norm{Av^+ + Bv^{0-}}.
  \end{equation}
  Since $A$ is a unipotent matrix, $\norm{A}\geq 1$.  We put
  \begin{equation} 
    \label{eq:kappa}
    \kappa=(1/3)\min\{1,\norm{B\inv}\inv\norm{A}\inv\}\leq (1/3)
    \min\{1,\norm{B\inv}\inv\}.
  \end{equation}
  Now to prove \eqref{eq:abs1} it is enough to consider the case when
  \begin{equation}
    \label{eq:vplus}
    \norm{v^+}\leq \kappa\|v\| \leq (1/3)\|v\|.  
  \end{equation}
  In particular, 
  \begin{equation}
    \label{eq:4}
    \norm{v^{0-}}\geq \norm{v}-\norm{v^+}\geq \norm{v}-(1/3)\norm{v} =
    (2/3)\norm{v}.
  \end{equation}
  Then by \eqref{eq:abs2}, \eqref{eq:kappa}, \eqref{eq:4}, and
  \eqref{eq:vplus},
  \[
\begin{array}{ll}
    \norm{(\boldsymbol{u}v)^{+0}}&\geq \norm{Bv^{0-}}-\norm{Av^+} \\
    &\geq \norm{B^{-1}}\inv \norm{v^{0-}}-\norm{A}\norm{v^+}\\
    &\geq \norm{B^{-1}}\inv \norm{v^{0-}} - \kappa\norm{A}\norm{v}\\
    &\geq (\norm{B^{-1}}\inv-(3/2)\kappa\norm{A})\norm{v^{0-}}\\
    &\geq (1/2)\norm{B^{-1}}\inv \norm{v^{0-}} \\
    &\geq (1/2)\norm{B^{-1}}\inv (2/3)\norm{v}\\
    &\geq \kappa\norm{v}.
  \end{array}
  \]
\end{proof}

\begin{cor} 
\label{cor:main1} 
Let $V$ be a finite dimensional normed linear space.  Consider a liner
representation of $G=\SO(n,1)$ on $V$, where $n\geq 2$.  Let
\begin{equation} 
\label{eq:V+-}
\begin{array}{ll}
V^+&=\{v\in V:a^{-k} v\stackrel{k\to\infty}{\longrightarrow}\infty,\
\forall a\in A^+\}\\
V^-&=\{v\in V:a^{k} v\stackrel{k\to\infty}{\longrightarrow}\infty,\
\forall a\in A^+\}\\
V^0&=\{v\in V:Av=v\}.
\end{array}
\end{equation}
Then given a compact set $F\subset N\setminus\{e\}$, there exists a
constant $\kappa>0$ such that for any $\boldsymbol{u}\in F$,
  \begin{equation} 
\label{eq:F-kappa}
    \max\{\norm{v^{+}},\norm{(\boldsymbol{u}v)^{+0}}\}
    \geq \kappa\norm{v}, \quad \forall v\in V.
  \end{equation}
  In particular, for any $a\in A^+$, and any $\boldsymbol{u}\in F$,
  \[
  \max\{\norm{av},\norm{a\boldsymbol{u}v}\}\geq \kappa\norm{v}, \quad
  \forall v\in V.
  \]
\end{cor}

\begin{proof}
  Given any $a\in A^+$ and $\boldsymbol{u}\in F$, there exist a
  continuous homomorphism of $\SL(2,\R)$ into $G$ such that $a$ is the
  image of $\psmat{\alpha & \\ & \alpha\inv}$ for some $\alpha>1$, and
  $\boldsymbol{u}$ is the image of $\psmat{1 & t \\0 & 1}$ for some
  $t\neq 0$.  We apply Lemma~\ref{lem:SL2} to obtain a constant
  $\kappa_1>0$ such that \eqref{eq:F-kappa} holds for
  $\boldsymbol{u}$.

  Now there exists a compact set $F_1\subset \cen(A)$ such that any
  $\boldsymbol{u}_1\in F$ is of the form $z\boldsymbol{u}z\inv$ for
  some $z\in F_1$.  Also there exists a constant $\kappa_2>0$ such that
  \[
  \kappa_2\norm{v}\leq \norm{zv}\leq \kappa_2\inv \norm{v}\quad
  \forall z\in F_1,\ \forall v\in V.
  \]
  Using this fact, we see that \eqref{eq:F-kappa} holds for any
  $\boldsymbol{u}_1\in F$ in place of $\boldsymbol{u}$ and
  $\kappa:=\kappa_2^2\kappa_1$.
\end{proof}

\subsection{Proof of Theorem~\ref{thm:return}}
Let $t_1,t_2\in I$ be such that
$\boldsymbol{u}:=u(\phi(t_2)-\phi(t_1))\inv\neq e$.
By Corollary~\ref{cor:main1} there exists $\kappa>0$ such that
\begin{equation}
  \label{eq:5}
  \sup{\norm{a_iv},\norm{a_i\boldsymbol{u}v}}\geq \kappa\norm{v},\quad \forall v\in V.
\end{equation}

Let a sequence $g_i\to g\in G$ be such that $\pi(g_i)=x_i$.  By
Proposition~\ref{prop:return} 
$\Gamma \Sigma \boldsymbol{p}$ is discrete in $V$. Therefore
\[
R_1:=\inf\{\norm{u(\phi(t_1)g_i\gamma\sigma
  \boldsymbol{p}}:\gamma\in\Gamma,\ \sigma\in\Sigma\}>0.
\]
For any $\gamma\in \Gamma$ and $\sigma\in \Sigma$, if we put
$v=u(\phi(t_1))g_i\gamma\sigma \boldsymbol{p}$ in \eqref{eq:5}, then
have
\begin{equation}
  \label{eq:1}
  \sup_{t_1,t_2}\{\norm{a_iu(\phi(t))g_i\gamma\sigma \boldsymbol{p}}\}\geq
  \kappa \norm{u(\phi(t_1))g_i\gamma\sigma \boldsymbol{p}}\geq \kappa R_1
\end{equation}

Now given $\epsilon>0$, and we obtain a compact set $\cK\subset\gmg$
such that the conclusion of Proposition~\ref{prop:return} holds for
$R=(1/2)\kappa R_1$.  Then by \eqref{eq:1}, for any $i\in\N$, the
possibility~(I) in the conclusion of Proposition~\ref{prop:return}
does not hold for $h_1=a_i$, $h_2=g_i$. Therefore the possibility~(II)
of Proposition~\ref{prop:return} must hold for all $i$.  Thus
Theorem~\ref{thm:return} follows.  \qed

We obtain the following immediate consequence of
Theorem~\ref{thm:return}:

\begin{cor} 
\label{cor:main-measure} 
After passing to a subsequence, $\mu_i\to\mu$ in the space of
probability measures on $\gmg$ with respect to the
weak$^\ast$-topology; that is,
  \[
  \lim_{i\to\infty}\int_{\gmg} f\,{\dd}\mu_i=\int_\gmg f\,{\dd}\mu,\quad
  \forall f\in\Cc(\gmg).
  \]
\end{cor}

\section{Invariance under a unipotent flow}

Let $I=[a,b]\subset\R$ with $a<b$.  Let $\phi:I\to \R^{n-1}$ be a
$\Ctwo$-curve such that $\dot\phi(t)\neq 0$ for all $t\in I$, where
$\dot\phi(t)$ denotes the tangent to the curve $\phi$ at $t$.  Fix
$w_0\in\R^{n-1}\setminus \{0\}$, and define
\[
W=\{u(tw_0):t\in\R\}.
\]
Consider the $\cen(A)$-action on $\R^{n-1}$ via the correspondence
$u(zv)=zu(v)z\inv$ for all $v\in\R^{n-1}$.  Then $\cen(A)=MA$ acts
transitively on $\R^{n-1}\setminus\{0\}$.  Therefore there exists a
continuous function $z:I\to \cen(A)$ such that
\begin{equation}
\label{eq:zw}
z(t)\dot\phi(t)=w_0, \quad \forall t\in I.
\end{equation}

Let a sequence $\{a_i\}_{i\in\N} \subset A^+$ be such that
$\alpha(a_i)\to\infty$ as $i\to\infty$. Let $x_i\to x$ a convergent
sequence in $G/\Gamma$. For each $i\in\N$, let $\lambda_i$ be the
probability measure on $\gmg$ such that
\begin{equation}
  \label{eq:lambda_i}
  \int_{\gmg} f d\lambda_i = \frac{1}{\abs{I}}\int_{t\in I}
  f(z(t)a_iu(\phi(t))x_i)\,{\dd}t,\quad \forall f\in \Cc(\gmg).
\end{equation}
Since $z(I)$ is compact, by Theorem~\ref{thm:return}, there exists a
probability measure $\lambda$ on $\gmg$ such that, after passing to a
subsequence, $\lambda_i\to \lambda$ as $i\to\infty$, in the space of
finite measures on $\gmg$ with respect to the weak$^\ast$-topology.

\begin{theo} 
  \label{thm:W-invariant}
  The measure $\lambda$ is $W$-invariant.
\end{theo}

\newcommand{\simm}[1]{\ \stackrel{#1}{\approx}\ }

\begin{proof}
  We will use the notation $\eta_1\simm{\epsilon} \eta_2$ to say
  $\abs{\eta_1-\eta_2}\leq \epsilon$.

  Let $f\in\Cc(\gmg)$ and $\epsilon>0$ be given.  Let $\Omega$ \nbd
  \begin{equation}
    \label{eq:Omega}
    f(\omega y)\simm{\epsilon} f(y),\quad 
    \forall\/ \omega\in\Omega^2\text{ and }\forall\/ y\in\gmg.  
  \end{equation}
 
  Let $t_0\in\R$. Let $t\in I=[a,b]$ and $i\in \N$. By
  \eqref{eq:zw}
\begin{equation}
  \label{eq:7}
  u(t_0w_0)z(t)a_i = z(t)a_iu(\alpha(a_i)\inv
  t_0z(t)\inv\cdot w_0)=z(t)a_iu(\xi_i\dot\phi(t)),
\end{equation}
where $\xi_i:=\alpha(a_i)\inv t_0$. Since $\phi$ is a
$\mathrm{C}^2$-map, 
\begin{equation}
  \label{eq:9}
  \phi(t+\xi_i)=\phi(t)+\xi_i\dot\phi(t)+\epsilon_i(t),
\end{equation}
where by Taylor's formula, there exists a
constant $M>0$ such that
\begin{equation}
  \label{eq:15}
\abs{\epsilon_i(t)}\leq M\abs{\xi_i}^2\leq
(M\abs{t_0}^2)\alpha(a_i)^{-2}, \quad\forall t\in[a,b].
\end{equation}

As $i\to\infty$, we have $\alpha(a_i)\to\infty$, and hence $\xi_i\to
0$ and $\alpha(a_i)\epsilon_i(t)\to 0$. Since $t\mapsto z(t)$ is
continuous, there exists $i_0\in\N$ such that for all $i\geq i_0$,
\begin{align}
  \label{eq:16}
  z(t+\xi_i)z(t)\inv \in \Omega\quad\text{and}\quad 
u(z(t)\cdot(\alpha(a_i)\epsilon_i(t)))\in\Omega.
\end{align}
Therefore 
\begin{align}
\label{eq:17}
\begin{array}{cl}
  &z(t+\xi_i)a_iu(\phi(t+\xi_i))\\
  =&(z(t+\xi_i)z(t)\inv)z(t)a_iu(\phi(t)+\xi_i\dot\phi(t)+\epsilon_i(t)),
  \quad \text{by \eqref{eq:9}}\\
  \in& 
  \Omega u(z(t)\cdot(\alpha(a_i)\epsilon_i(t)))
  z(t)a_i(u(\phi(t)+\xi_i\dot\phi(t)))\\
  \subset& \Omega^2 z(t)a_iu(\xi_i\dot\phi(t))u(\phi(t)),\quad \text{by
    \eqref{eq:16}}\\
  \subset &\Omega^2 u(t_0w_0)z(t)a_iu(\phi(t)), \quad \text{by \eqref{eq:7}}.  
\end{array}
\end{align}
Therefore by \eqref{eq:Omega} 
\begin{equation}
  \label{eq:19}
  f(z(t+\xi_i)a_iu(\phi(t+\xi_i)x_i)
  \simm{\epsilon} f(u(t_0w_0)z(t)a_iu(\phi(t))x_i).
\end{equation}
Hence
\begin{align}
  \label{eq:18}
\begin{array}{cl}
&\int_a^b f(z(t)a_iu(\phi(t))x_i)\,{\dd}t\\
\simm{\xi_i\sup \abs{f}} &\int_a^{b-\xi_i}
f(z(t+\xi_i)a_iu(\phi(t+\xi_i))x_i)\,{\dd}t\\
\simm{\epsilon\abs{I}}&\int_a^{b-\xi_i}
f(u(t_0w_0)z(t)a_iu(\phi(t+\xi_i))x_i)\,{\dd}t,\quad \text{by \eqref{eq:19}}\\
\simm{\xi_i\sup\abs{f}} &\int_a^b f(u(s_0w_0)z(t)a_iu(\phi(t))x_i)\,{\dd}t.
\end{array}
\end{align}

Therefore, since $\epsilon>0$ is chosen arbitrarily, and $\xi_i\to 0$
as $i\to0$,
\begin{equation}
  \label{eq:20}
  \int_\gmg f(u(s_0w_0)y)\,{\dd}\lambda(y)=\int_\gmg f(y)\,{\dd}\lambda(y).
\end{equation}
\end{proof}

The above simplification of the original proof given in
arXiv:0708.4093v1 is based on referee's suggestions.

\section{Dynamical behaviour of translated trajectories \\ near
  singular sets}
\label{sec:linearization}

Let the notation be as in the previous section.  We will further assume
that $\phi:I\to \R^{n-1}$ is an analytic function.   In this case we
will further observe that the function $z:I\to \cen(A)$ such that
$z(t)\dot\phi(t)=w_0$ for all $t\in I$ is also an analytic
function.  Given a convergent sequence $x_i\to x$ in $\gmg$, we obtain
a sequence of measures $\{\lambda_i:i\in\N\}$ on $\gmg$ as defined by
\eqref{eq:lambda_i}.  Due to Theorem~\ref{thm:return}, by passing to a
subsequence we will assume that $\lambda_i\to\lambda$ as $i\to\infty$,
where $\lambda$ is a probability measure on $\gmg$.   By
Theorem~\ref{thm:W-invariant}, $\lambda$ is invariant under the
action of the one-parameter subgroup $W=\{u(sw_0):s\in\R\}$.  We would
like to describe the measure $\lambda$ using the description of
ergodic invariant measures for unipotent flows on homogeneous spaces
due to Ratner~\cite{R:measure}.  We begin with some notation.

Let $\sH$ denote the collection of analytic subgroups $H$ of $G$ such
that $H\cap\Gamma$ is a lattice in $H$, and a unipotent one-parameter
subgroup of $H$ acts ergodically with respect to the $H$-invariant
probability measure on $H/H\cap\Gamma$.  Then $\sH$ is a countable
collection \cite{Shah:uniform,R:measure}.

For $H\in\sH$, define
\begin{align*}
  N(H,W)&=\{g\in G: g\inv Wg\subset H\} \qquad\textrm{and}\\
  S(H,W)&=\tcup{\substack{F\in\sH\\F\subsetneq H}}N(F,W).
\end{align*}
Then by \cite[Lemma 2.4]{Moz+Shah:limit}
\begin{equation}
  \label{eq:piNS}
  \pi(N(H,W)\setminus S(H,W))=\pi(N(H,W))\setminus \pi(S(H,W)).
\end{equation}

Then by Ratner's theorem~\cite{R:measure}, as explained in
\cite[Theorem 2.2]{Moz+Shah:limit}:  
\begin{theo}[Ratner]
  \label{thm:Ratner}
  Given the $W$-invariant probability measure $\lambda$ on $\gmg$,
  there exists $H\in\sH$ such that
  \begin{equation} 
    \label{eq:lambda-H} 
    \lambda(\pi(N(H,W))>0 \quad\textrm{and} \quad \lambda(\pi(S(H,W))=0.
\end{equation}
Moreover almost every $W$-ergodic component of $\lambda$ on
$\pi(N(H,W))$ is a measure of the form $g\mu_H$, where $g\in
N(H,W)\setminus S(H,W)$, $\mu_H$ is a finite $H$-invariant measure on
$\pi(H)\cong H/H\cap\Gamma$, and $g\mu_H(E):=\mu(g\inv E)$ for all
Borel sets $E\subset \gmg$.
\end{theo}

For $d=\dim H$, let $V=\wedge^d \la{G}$, and consider the
$\wedge^d\Ad$-action of $G$ on $V$.  Fix
$p_H=\wedge^d\la{H}\setminus\{0\}$.

We recall some facts from \cite[\S3]{Moz+Shah:limit}: For any $g\in
\nor(H)$, $gp_H=\det(\Ad g|_{\la{H}})p_H$. Hence the stabilizer of
$p_H$ in $G$ equals
\[
\norone(H):=\{g\in \nor(H):\det((\Ad g)|_{\la{H}})=1\}.
\]
Since, $(\Gamma\cap \nor(H))\pi(H)=\pi(H)$, we have $(\Gamma\cap
\nor(H))p_H=p_H$ or $(\Gamma\cap \nor(H))p_H=\{p_H,-p_H\}$.  In the
former case we put $\bar V=V$ and in the later case we put $\bar
V=V/\{\pm 1\}$.  For any $v\in V$, we denote by $\bar v$ the image of
$v$ in $\bar V$, and define the action of $g\in G$ by $g\cdot\bar
v:=\bar{gv}$. We define $\eta(g)=g\bar p_H$ for all $g\in G$.

\begin{prop}[\cite{Dani+Mar:limit}] 
  \label{prop:eta-Gamma}
  $\eta(\Gamma)=\Gamma \cdot \bar p_H$ is a discrete subset of $\bar
  V$.  \qed
\end{prop}

Let $\cA=\{\bar v\in \bar V: v\wedge w_0=0\in
\wedge^{d+1}\la{g}\}$. Then $\cA$ is the image of a linear subspace of
$V$. We observe that
\begin{equation}
  \label{eq:A}
  N(H,W)=\eta\inv(\cA).
\end{equation}

Given any compact set $\cD\subset \cA$, we define
\begin{displaymath}
  S(\cD)=\{g\in \eta\inv(\cD):\eta(g\gamma)\in \cD \quad \textrm{for
    some } \gamma\in \Gamma\setminus \nor(H)\}.  
\end{displaymath}

\begin{prop}[{\cite[Proposition 3.2]{Moz+Shah:limit}}]
  \label{prop:inject}
  (1) $S(\cD)\subset S(H,W)$ and $\pi(S(\cD))$ is closed in $\gmg$.
  (2) For any compact set $\cK\subset \gmg \setminus \pi(S(\cD))$,
  there exists a neighbourhood $\Phi$ of $\cD$ in $\bar V$ such that
  for any $g\in G$ and $\gamma_1,\gamma_2\in\Gamma$:
  \begin{equation} 
    \label{eq:inject} 
    \textrm{if $\pi(g)\in \cK$ and
      $\{\eta(g\gamma_1),\eta(g\gamma_2)\}\subset \cl{\Phi}$, then
      $\eta(\gamma_1)=\eta(\gamma_2)$,}
  \end{equation}   
where $\cl{\Phi}$ denotes the closure of $\Phi$ in $\bar V$. 
\end{prop}

\subsection{$(C,\alpha)$-good family} 
\label{subsec:CAlpha-2} 
Let $\sF$ denote the $\R$-span of all the coordinate functions of
the maps $t\mapsto (\wedge^d\Ad)(z(t)u(\phi(t)))$ from $I$ to
$\GL(V)$.  As explained in \S\ref{subsec:CAlpha-1}, by
\cite[Proposition~3.4]{Klein+Mar:Annals98}, the family $\sF$ is
`$(C,\alpha)$-good' for some $C>0$ and $\alpha>0$; that is, for any
subinterval $J\subset I$, $\xi\in\sF$, and $r>0$,
\begin{equation} 
  \label{eq:CAlpha2} 
  \abs{\{t\in J: \abs{\xi(t)}<r\}} <
  C\left(\frac{r}{\sup_{t\in J}|\xi(t)|}\right)^\alpha\abs{J}.
\end{equation}

\begin{prop}[Cf.~\cite{Dani+Mar:limit}]
  \label{prop:Phi-Psi}
  Given a compact set $\cC\subset \cA$ and $\epsilon>0$, there exists
  a compact set $\cD\subset\cA$ containing $\cC$ such that given any
  neighbourhood $\Phi$ of $\cD$ in $\bar V$, there exists a
  neighbourhood $\Psi$ of $\cC$ in $\bar V$ contained in $\Phi$ such
  that for any $h\in G$, any $v\in \bar V$ and any interval $J\subset
  I$, one of the following holds:
  \begin{enumerate}
  \item[I)]$hz(t)u(\phi(t))v\in \Phi$ for all $t\in J$.
  \item[II)]
    $\abs{\{t\in J: hz(t)u(\phi(t))v\in\Psi\}} \leq \epsilon 
      \abs{\{t\in J:hz(t)u(\phi(t))v\in\Phi\}}$.
  \end{enumerate}
\end{prop}

\begin{proof}
  The argument in the proof of \cite[Proposition~4.2]{Dani+Mar:limit}
  goes through with straightforward changes.  Since $\cA$ is the image
  of a linear subspace of $V$, one can describe the neighbourhoods of
  subsets of $\cA$ in $\bar V$ via linear functionals.  Further, one uses
  the property \eqref{eq:CAlpha2} of the functions in $\sF$ instead of
  \cite[Lemma~4.1]{Dani+Mar:limit} in the proof.
\end{proof}

\subsection{Linear presentation of dynamics in thin neighbourhoods of
  singular sets}
\label{subsec:nbd}

Now let $C$ be any compact subset of $N(H,W)\setminus S(H,W)$.  Let an
$\epsilon>0$ be given.  We apply Proposition~\ref{prop:Phi-Psi} to
$\cC:=\eta(C)\subset\cA$, and obtain a compact set $\cD\subset \cA$.
By \eqref{eq:piNS}, since $\pi(C)$ is a compact subset
$\gmg\setminus\pi(\cS(\cD))$. We choose a compact set $\cK\subset
\gmg\setminus\pi(\cS(\cD))$ such that $\pi(C)$ is contained in the
interior of $\cK$.  Then we take any neighbourhood $\Phi_1$ of $D$ in
$\bar V$.  By Proposition~\ref{prop:inject}, there exists an open
neighbourhood $\Phi$ of $\cD$ contained in $\Phi_1$ such that
property \eqref{eq:inject} holds.  Now we obtain a neighbourhood $\Psi$ of
$\cC$ in $\bar V$ such that the conclusion of
Proposition~\ref{prop:Phi-Psi} holds.  Let
\begin{equation}
  \label{eq:22}
  \cO:=\pi(\eta\inv(\Psi))\cap\cK.
\end{equation}

Then $\cO$ is a neighbourhood of $\pi(C)$ in $\gmg$.

\begin{prop}[Cf.~\cite{Moz+Shah:limit}]
  \label{prop:main3}
  For any $h_1,h_2\in G$, and for any subinterval $J\subset I$, one of
  the following holds:
  \begin{enumerate}
  \item[a)] There exists $\gamma\in\Gamma$ such that
    $\eta(h_1z(t)u(\phi(t))h_2\gamma)\in\Phi$, $\forall t\in J$.
  \item[b)] $\abs{\{t\in J:\pi(h_1z(t)u(\phi(t))h_2)\in\cO\}}\leq
    (2\epsilon)|J|$.
  \end{enumerate}
\end{prop}

\begin{proof}
  Suppose that the possibility~(a) does not hold for the given $J$.

  Let $\psi(t):=h_1z(t)u(\phi(t))h_2$ for all $t\in I$.  Let
  \begin{equation}
    \label{eq:JO}
  J^\ast=\{t\in J:\pi(\psi(t))\in\cO\}.
\end{equation}
Take any $t\in J^\ast$.  By the choice of $\Phi$ with
property~\eqref{eq:inject}, there exists a unique $v_t\in\eta(\Gamma)$
such that $\psi(t)v_t\in\cl{\Phi}$; and hence due to \eqref{eq:22} and
\eqref{eq:JO} we have $\psi(t)v_t\in \Psi$.  Let $J(t)$ be the
largest subinterval of $J$ containing $t$ such that
\begin{equation}
\label{eq:Jt}
\psi(s)v_t\in \cl{\Phi},\quad \forall s\in J(t).  
\end{equation}
By \eqref{eq:inject} and \eqref{eq:Jt}, we have 
\begin{equation}
  \label{eq:21}  
  v_s=v_t, \text{ and hence } \psi(s)v_t=\psi(s)v_s\in\Psi, \quad \forall s\in J^\ast\cap J(t). 
\end{equation}
 
Since the possibility~(a) does not hold for $J$, by our choice of
$J(t)$, we have that $J(t)$ contains one of its end-points, say $s_e$,
and $\psi(s_e)v_t\not\in \Phi$. Thus $\psi(J(t))v_t\not\subset \Phi$.
Therefore by Proposition~\ref{prop:Phi-Psi}, in view of \eqref{eq:22}
and \eqref{eq:21}, we deduce that
  \begin{equation}
    \label{eq:compare}
    \abs{J^\ast\cap J(t)}\leq \epsilon\abs{J(t)}.
  \end{equation}
 
  Due to \eqref{eq:21}, $J(s)=J(t)$ for all $s\in J^\ast\cap
  J(t)$. Therefore there exists a countable set $\cJ^\ast\subset
  J^\ast$ such that
\begin{equation} 
\label{eq:Jast} 
J^\ast\subset \tcup{t\in \cJ^\ast} J(t),
\end{equation}
and if $t_1\neq t_2$ in $\cJ^\ast$ then $t_1\not\in J(t_2)$.  

In particular, if $t_1<t_2$ in $\cJ^\ast$ then $J(t_1)\cap J(t_2)\subset
(t_1,t_2)$.  Therefore if $t_1<t_2<t_3$ in $\cJ^\ast$, then 
\[
J(t_1)\cap J(t_2)\cap J(t_3) = \emptyset.
\]
 Hence 
 \begin{equation} 
\label{eq:double}
 \sum_{t\in\cJ^\ast} \abs{J(t)} \leq 2 \abs{\tcup{t\in\cJ^\ast}
   J(t)}.
\end{equation}

 Now by \eqref{eq:compare}, \eqref{eq:Jast} and \eqref{eq:double},
  \[
  \abs{\cJ^\ast} \leq \epsilon \sum_{t\in\cJ^\ast} \abs{J(t)} \leq
  (2\epsilon)\abs{J}.
  \]
\end{proof}

\subsection{Algebraic consequences of positive limit measure on
  singular sets}

\label{subsec:conseq}
Let $\{a_i\}\subset A$ and
$x_i\stackrel{i\to\infty}{\longrightarrow}x$ be the sequences
which are involved in the definition of $\lambda_i$,
see~\eqref{eq:lambda_i}.

Let $V^+=\{v\in V:a_i\inv v\stackrel{i\to\infty}{\longrightarrow}
0\}$, $V^0=\{v\in V:Av=v\}$, and $V^-=\{v\in V: a_i
v\stackrel{i\to\infty}{\longrightarrow}0\}$. Then $V=V^+\oplus
V^0\oplus V^-$. Let $q^+:V\to V^+$ and $q^{+0}:V\to V^++V^0$ denote
the associated projections. Let $\bar V^+$, $\bar V^0$ and $\bar V^-$
denote the projections of $V^+$, $V^0$ and $V^-$, on $\bar V$
respectively. The sets $V^\pm$ do not change if we pass to a
subsequence of $\{a_i\}$.

We recall that after passing to a subsequence, $\lambda_i\to \lambda$
in the space of probability measures on $G/\Gamma$, and by
Theorem~\ref{thm:W-invariant} and Theorem~\ref{thm:Ratner} there
exists $H\in\sH$ such that
\begin{equation}
  \label{eq:13}
\lambda(\pi(N(H,W)\setminus S(H,W)))>0.
  \end{equation}

  The goal of this section is to analyze this condition using
  Proposition~\ref{prop:main3} and Corollary~\ref{cor:main1} to obtain its
  following algebraic consequence.

\begin{prop} 
\label{prop:V0-}
Let $g\in G$ be such that $\pi(g)=x$. Then there exists
$\gamma\in\Gamma$ such that
\begin{equation} 
  \label{eq:V0minus}
  \eta(z(t)u(\phi(t))g\gamma) \subset \bar V^0+\bar V^-, \quad \forall t\in I.
\end{equation}
\end{prop}

\begin{proof}
By \eqref{eq:13} there exists a compact
set $C\subset N(H,W)\setminus S(H,W)$ such that
$\lambda(\pi(C))>c_0>0$ for some constant $c_0>0$.  We fix
$0<\epsilon<c_0/2$, and obtain the compact sets $\cD\subset\cA$ and
$\cK\subset \gmg$ as in \S\ref{subsec:nbd}. Next we choose any
neighbourhood $\Phi_1$ of $\cD$ in $\bar V$, and obtain a
neighbourhood $\Psi$ of $\eta(C)$ as in \S\ref{subsec:nbd}. Let
$i_1\in\N$ be such that if we put $\cO:=\pi(\eta\inv(\Psi))\cap\cK$,
then
\begin{equation} 
  \label{eq:cO}
\lambda_i(\cO)>c_0 \quad \textrm{for all $i\geq i_1$}.
\end{equation}

Since $x_i\xrightarrow{i\to\infty} x$ and $\pi(g)=x$, there exists a
convergent sequence $g_i\stackrel{i\to\infty}{\longrightarrow} g$ in
$G$ such that $\pi(g_i)=x_i$ for all $i\in\N$.  By
\eqref{eq:cO} and \eqref{eq:lambda_i}, since $z(t)\in\cen(A)$,
\begin{equation}
  \label{eq:Olambda}
  \abs{\{t\in I: \pi(a_iz(t)u(\phi(t))g_i)\in\cO\}}> c_0\abs{I}, \quad
  \forall i\geq i_1.
\end{equation}

We apply Proposition~\ref{prop:main3} for $h_1=a_i$, $h_2=g_i$, and
$J=I$.  Then since $c_0>2\epsilon$, by \eqref{eq:Olambda},
the possibility~(b) in the conclusion of Proposition~\ref{prop:main3}
does not hold, and hence possibility~(a) of the proposition must hold;
that is, there exists $\gamma_i\in\Gamma$ such that
\begin{equation}
  \label{eq:gamma_i}
  \eta(a_iz(t)u(\phi(t))g_i\gamma_i)\in\Phi\subset\Phi_1, \quad\forall
  i\geq i_1.
\end{equation}

We choose a decreasing sequence of neighbourhoods $\Phi_k$ of $\cD$ in
$\bar V$ be such that $\bigcap_{k\in\N} \Phi_k=\cD$, and apply the
above argument for each $\Phi_k$ in place of $\Phi_1$. We then obtain
sequences $i_k\stackrel{k\to\infty}{\longrightarrow}\infty$ in $\N$
and $\{\gamma_k\}$ in $\Gamma$ such that
\[
a_{i_k}z(t)u(\phi(t))\eta(g_{i_k}\gamma_k)\in\Phi_k, \quad \forall
t\in I, \ \forall k\in\N.
\]

Since $\{z(t):t\in I\}$ is contained in a compact set, there exists
$R>0$ such that $z(I)\inv \Phi_1$ is contained in $B(R)$, the ball of
radius $R$ centered at $0$ in $\bar V$.  Thus
\begin{equation}
\label{eq:R1}
\norm{a_{i_k}u(\phi(t))\eta(g_{i_k}\gamma_{k})}\leq R, \qquad \forall
t\in I.
\end{equation}

Fix any $t_1\in I$. Since $\phi$ is a nonconstant function, by
Corollary~\ref{cor:main1}, there exists a constant $\kappa>0$ such
that
\begin{equation}
\label{eq:1q0plus}
\sup_{t\in I} \norm{q^{+0}(u(\phi(t))v)}\geq
\kappa\norm{u(\phi(t_1))v},\quad \forall v\in V.
\end{equation}

For all $v\in V$, define $\norm{\bar v}:=\norm{v}$, and let $q^{+0}(\bar
v)$ and $q^+(\bar v)$ denote the images of $q^{+0}(v)$ and $q^+(v)$ in
$\bar V$, respectively.  Let $t\in I$. Then
\begin{align*}
&\norm{q^{+0}(u(\phi(t)\eta(g_{i_k}\gamma_{k})))}\\
&\leq \norm{a_{i_k}q^{+0}(u(\phi(t)\eta(g_{i_k}\gamma_{k})))}\\
&\leq \norm{a_{i_k}u(\phi(t))\eta(g_{i_k}\gamma_{k})}\\
&\leq R, \qquad \text{(by \eqref{eq:R1})}.
\end{align*}
Therefore by \eqref{eq:1q0plus},
\begin{equation}
\label{eq:uinv}
\norm{\eta(g_{i_k}\gamma_{k})}\leq \kappa \inv \norm{u(\phi(t_1))\inv}.
\end{equation}
Since $\eta(\Gamma)$ is discrete and
$g_{i_k}\stackrel{k\to\infty}{\longrightarrow} g$, due to
\eqref{eq:uinv}, the set $\{\eta(\gamma_k):k\in\N\}$ is finite.
Therefore by passing to a subsequence, there exists $\gamma\in\Gamma$
such that $\eta(\gamma_k)=\eta(\gamma)$ for all $k\in\N$ and hence
\begin{equation} 
\label{eq:gamma2} 
\eta(a_{i_k}z(t)u(\phi(t))g_{i_k}\gamma)\in\Phi_k,\ \forall k\in\N.
\end{equation}
For each $k\in\N$, if $w_k^+=q^+(\eta(z(t)u(\phi(t))g_{i_k}\gamma))\in
\bar{V}^+$, then by \eqref{eq:gamma2} we have $\limsup_{k\to\infty}
\norm{a_{i_k}w_k^+}<\infty$. Since
$\alpha(a_{i_k})\stackrel{k\to\infty}{\longrightarrow}\infty$, we
conclude that $w_k^+\stackrel{k\to\infty}{\longrightarrow}0$. Since
$g_{i_k}\stackrel{k\to\infty}{\longrightarrow} g$, we have
\[
q^+(\eta(z(t)u(\phi(t)g\gamma)))=\lim_{k\to\infty}
q^+(\eta(z(t)u(\phi(t)g_{i_k}\gamma)))=\lim_{k\to\infty}w_k^+=0.
\]
Therefore \eqref{eq:V0minus} follows.
\end{proof}

In order to derive group theoretic consequences of condition
\eqref{eq:V0minus} we will need the following observation.

\begin{lemma} 
  \label{lema:1} 
  If $\eta(h_1),\eta(h_2)\in \bar V^0\cap\cA$ for some $h_1,h_2\in G$, then
  \[
  A\subset h_1Hh_1\inv,\quad \nor(H)^0\subset \cen(H)H,
  \quad\text{and} \quad h_2\in (M\cap \cen(W))h_1H.
  \]

  Also $h_1Hh_1\inv$ contains a cocompact normal subgroup containing
  $A$ which is conjugate to $\SO(m,1)$ in $G$, where $2\leq m\leq n$.
\end{lemma}

\begin{proof}
  Since $A$ is connected, $A\subset h_j\norone(H)h_j\inv$,
  $j=1,2$. Let $H_1$ denote the Zariski closure of $H$ in $G$. Then
  $A\subset h_j\norone(H_1)h_j\inv$, $j=1,2$. If $U_1$ denotes the
  unipotent radical of $H_1$, then
\begin{equation} \label{eq:AU}
A\subset h_j\norone(U_1)h_j\inv, \quad j=1,2.
\end{equation}
Any nontrivial unipotent element of $G$ is contained in a unique
maximal unipotent subgroup of $G$. Therefore either
$h_1U_1h_1\inv\subset N$ or $h_1U_1h_1\inv \subset N^-$ (recall
\eqref{eq:N}). Therefore in view of \eqref{eq:AU}, we conclude that
$U_1=\{e\}$.  Therefore $\nor(H_1)$ is reductive, and hence
$\nor(H)=\nor(H_1)$ is reductive. Therefore
$\nor(H)^0=\cen(H)^0H\subset\norone(H)$. Since $\nor(H)$ contains
nontrivial unipotent elements, its maximal semisimple factor contains
an $\R$-diagonalizable subgroup. Therefore, since $G$ is of $\R$-rank
one, $\cen(H)$ is compact. In particular, $A\subset h_jHh_j\inv$,
$j=1,2$.

Being a reductive subgroup of $\SO(n,1)$ containing a nontrivial
unipotent element, each $h_jHh_j\inv$ is of the form
$vM_1\SO(m,1)v\inv$ for some $v\in N$ (see \eqref{eq:N}), $2\leq m\leq
n$ and $M_1\subset \cen(\SO(m,1))\subset K$.  Note that if $g_i\inv AW
g_i\subset \SO(m,1)$ for some $g_1,g_2\in G$, then there exists $g'\in
\SO(m,1)$, such that $(g')\inv(g_1\inv aw g_1)g' = g_2\inv aw g_2$ for
all $a\in A$ and $w\in W$. Therefore $g_2\in \cen(AW)g_1g'\subset
(M\cap\cen(W))g_1\SO(m,1)$. Thus in our situation, we conclude that
have $h_2\in (M\cap\cen(W))h_1H$.
\end{proof}

\subsection{Algebraic description of $\lambda$}
\begin{prop}
\label{prop:lambda}
There exist an analytic map $\tilde\xi:I\to \cen(W)\cap M$ and an
element $h_1\in G$ such that $AWh_1\subset h_1H$ and the following
holds: For any $f\in\Cc(G/\Gamma)$, we have
  \begin{equation} 
    \label{eq:lambda-xi}
    \int_{G/\Gamma} f\,{\dd}\lambda = \frac{1}{\abs{I}}\int_{t\in I}\left(
      \int_{y\in\pi(H)}
      f(\tilde\xi(t)h_1y)\,{\dd}\mu_H(y)\right){\dd}t.
  \end{equation}
In particular, $\lambda$ is $AW$-invariant. Moreover
\begin{equation}
  \label{eq:n-}
u(\phi(I))\subset N\cap (P^-h_1H(g\gamma)\inv)=N\cap (N^-Mh_1H(g\gamma)\inv).
\end{equation}
\end{prop}

\begin{proof}
  Let the notation be as in the previous section. We start with a
  construction.  Since $H\cap \Gamma$ is a lattice in $H$ and
  $\nor(H)/H$ is compact, we have that $(\nor(H)\cap\Gamma)/(H\cap
  \Gamma)$ is finite.  Since $\eta(\Gamma)$ is discrete, the map
  $\rho:G/(H\cap \Gamma)\to (\gmg)\times \bar V$, defined by
\[
\rho(h(H\cap \Gamma))=(\pi(h),\eta(h)),\quad \forall h\in G,
\]
is a continuous proper map.

By Proposition~\ref{prop:V0-} there exists $\gamma\in\Gamma$ such that
\[
  \eta(z(t)u(\phi(t))g\gamma) \subset \bar V^0+\bar V^-, \quad \forall t\in I.
\]
We put $\xi(t):= q^0(z(t)\eta(u(\phi(t))g\gamma))\in \bar V^0$ for all
$t\in I$.  Then $\xi:I\to \bar{V^0}$ is an analytic function, and
\begin{equation} 
\label{eq:converge-xi}
  \xi(t) = 
  \lim_{\substack{\alpha(a)\to\infty\\a\in A^+}}\eta(az(t)u(\phi(t))g\gamma)\in\cD.
\end{equation}
Thus 
\begin{equation}
  \label{eq:6}
  \xi(I)\subset \cA\cap \bar V^0.
\end{equation}

For each $i\in\N$, we define a probability
measure $\tilde\lambda_i$ on $G/(H\cap\Gamma)$ such that for any
$\tilde f\in\Cc(G/(H\cap \Gamma))$,
  \[
  \int_{G/(H\cap \Gamma)} \tilde
  f\,{\dd}\tilde\lambda_i=\frac{1}{\abs{I}}\int_I \tilde
  f(a_iz(t)u(\phi(t))g_i\gamma(H\cap \Gamma))\,{\dd}t.
  \]
  Then $\tilde\lambda_i$ projects onto $\lambda_i$ under the quotient
  map $\rho_1:G/(H\cap \Gamma)\to G/\Gamma$. Let
  $\rho_2:G/(H\cap\Gamma)\to \bar V$ be the map defined by
  $\rho_2(h(H\cap\Gamma))=\eta(h)$ for all $h\in G$. Then due to
  \eqref{eq:converge-xi}, the projected measures
  $(\rho_2)_\ast(\tilde\lambda_i)$ on $\bar V$ converge to a
  probability measure, say $\nu$, supported on $\xi(I)$ as
  $i\to\infty$. Now since $\rho$ is a proper map, we conclude that,
  after passing to a subsequence, as $i\to\infty$, $\tilde\lambda_i$
  converges to a probability measure $\tilde \lambda$ on $G/(H\cap
  \Gamma)$ such that
  \begin{equation}
    \label{eq:rho12}
    (\rho_1)_\ast(\tilde\lambda)=\lambda \quad \textrm{and} \quad
    (\rho_2)_\ast(\tilde\lambda)=\nu. 
  \end{equation}

  Therefore by \eqref{eq:6} and Lemma~\ref{lema:1}, there exists
  $h_1\in N(H,W)$ such that
\begin{equation}
    \label{eq:xi-NHW}
    \xi(t)\in (\cen(W)\cap M)\eta(h_1),\quad \forall t\in I.
  \end{equation}
Hence
  \begin{equation}
    \label{eq:tilde-lambda} 
    \supp(\tilde\lambda)\subset (\cen(W)\cap M)h_1H/(H\cap \Gamma).
  \end{equation}
    
  Since $G$ is an algebraic group acting linearly on $V$, the orbit
  $\eta(G)$ is open in its closure, and hence locally compact in the
  relative topology. Now $\nor(H)^0\subset \cen(H)H$ and
  $\cen(H)\subset h_1\inv M h_1$ is compact. In particular, the
  quotient map $G/H\to \eta(G)$ given by $hH\mapsto\eta(h)$, for all
  $h\in G$, is a proper map with respect to the relative topology on
  $\eta(G)\subset \bar V$. Therefore due to \eqref{eq:xi-NHW}, there
  exists an analytic map $\tilde\xi:I\to \cen(W)\cap M$ such that
  \begin{equation}
    \label{eq:tilde-xi}
    \lim_{\substack{\alpha(a)\to\infty\\ a\in A^+}}
    az(t)u(\phi(t))g\gamma H=\tilde{\xi}(t)h_1H,\qquad \forall t\in I.
  \end{equation}

  By \eqref{eq:rho12} and
  \eqref{eq:tilde-lambda} $\lambda$ is concentrated on
  $\pi(N(H,W))$. Then almost every normalized $W$-ergodic component of
  $\lambda$ is of the form $hh_1\mu_H$ for some $h\in \cen(W)\cap M$,
  where $\mu_H$ is the $H$-invariant probability measures on
  $\pi(H)$. Therefore, since $hh_1\mu_H$ is $A$-invariant for each
  $h\in M$, we conclude that $\lambda$ is $A$-invariant. 

  Let $\tilde\eta :G/(H\cap\Gamma)\to G/H$ be the quotient map.  Let
  $\bar\lambda=\tilde\eta_\ast(\tilde\lambda)$ be the projection of
  $\tilde\lambda$ on $G/H$. Then for any $\tilde
  f\in\Cc(G/H\cap\Gamma)$,
  \begin{equation} 
    \label{eq:bar-lambda} 
    \int_{G/(H\cap\Gamma)} \tilde f\,{\dd}\tilde\lambda 
    = \int_{G/H} \left(\int_{y\in\pi(H)}
      \tilde{f}(hh_1y)\,\dd\mu_H(y)\right) \,{\dd}\bar\lambda(hh_1H).
  \end{equation}
  By \eqref{eq:tilde-xi}, $\bar\lambda$ is the projection of the
  normalized Lebesgue measure of $I$ onto $\tilde\xi(I)h_1H/H$. Thus
  we obtain a complete description of the measure $\lambda$ as in
  \eqref{eq:lambda-xi}.

Since $h_1Hh_1\inv$ is a reductive subgroup of $G$ containing $A$, if
there exist $h'\in G$ and $h\in M$ such that
  \[
  \lim_{\substack{\alpha(a)\to\infty \\ a\in A^+}}ah'(h_1H)=h(h_1H),
  \]
  then $h'\in N^-h$.  Hence by \eqref{eq:tilde-xi}, there exists a
  continuous map $n^-(t):I\to N^-$ such that
  \begin{equation*}
  z(t)u(\phi(t))g\gamma\in n^-(t)\tilde\xi(t)(h_1H), \quad\forall t\in I.   
  \end{equation*}
  Therefore, since $Ah_1H=h_1H$, we obtain \eqref{eq:n-}.
\end{proof}

\section{Proofs of results stated in the Introduction}

In order to describe the limiting distributions for the sequence of
measures $\mu_i$ as defined in \eqref{eq:mui} using
Proposition~\ref{prop:lambda}, we make the following observation.

\begin{lemma}
  \label{lem:ZA-translate}
  Let $\{\theta_i:I=[a,b]\to G/\Gamma\}_{i\in\N}$ be sequence of
  continuous curves, and $\{a_i\}_\N\subset A^+$ be a sequence such
  that $\alpha(a_i)\stackrel{i\to\infty}{\longrightarrow}\infty$. Let
  $E\subset I$ be a finite set, and suppose that for each $t\in
  I\setminus E$ there exists a probability measure $\lambda_t$ on
  $G/\Gamma$ such that the map $t\mapsto \lambda_t$ is continuous on
  $I\setminus E$ with respect to the weak$^\ast$-topology on the space
  of probability measures on $G/\Gamma$, and for every closed interval
  $J\subset I\setminus E$ with nonempty interior, we have
  \begin{equation*}
    \lim_{i\to\infty} \frac{1}{\abs{J}}\int_J
    f(a_i\theta_i(t))\,{\dd}t = \int_{t\in J} \left(\int_\gmg
      f\,{\dd}\lambda_t \right){\dd}t.
  \end{equation*}
  Let $z_i\stackrel{i\to\infty}{\longrightarrow} z$ from $I$ to $P^-$,
  and and $w_i\stackrel{i\to\infty}{\longrightarrow} w$ from $I$ to
  $G$ be uniformly convergent sequences of continuous functions. Then
  \begin{equation*}
    \lim_{i\to\infty} \int_If(w_i(t)a_iz_i(t)\theta_i(t))\,{\dd}t 
    = \int_I\left(\int_\gmg
      f(w(t)\zeta(t)y)\,{\dd}\lambda_t(y)\right){\dd}t ,
  \end{equation*}
  where $\zeta(t)\in \cen(A)$ is such that $z(t)\in N^-\zeta(t)$ for
  all $t\in I$.

  In particular if $\lambda_t=\mu_G$ for all $t\in I$, then
  \begin{equation}
    \label{eq:12}
    \lim_{i\to\infty}
    \frac{1}{\abs{I}}\int_If(w_i(t)a_iz_i(t)\theta_i(t))\,{\dd}t 
    = \int_\gmg f\,{\dd}\mu_G.
  \end{equation}
  \qed
\end{lemma}

\begin{proof}[\bf Proof of Theorem~\ref{thm:psi}]
  Since $\phi$ is analytic and nonconstant, the set $E:=\{t\in
  I:\dot\phi(t)=0\}$ is finite. Let $J\subset I\setminus E$ be any
  closed interval with nonempty interior. Let $g_1\in\pi\inv(x)$. We
  consider the discussion of \S\ref{sec:linearization} for $J$ in
  place of $I$ and $g_i=g_1$ for all $i$. Then there exist a reductive
  subgroup $H\in\sH$ and $h_1\in G$ and $g\in g_1\Gamma$ such that
  $AW\subset h_1Hh_1\inv$ and by \eqref{eq:n-},
\begin{equation} 
  \label{eq:u-J}
  u(\phi(J))\subset N\cap (N^-Mh_1Hg\inv).
\end{equation}
Without loss of generality we may assume that $H$ is a smallest
dimensional subgroup of $\sH$ such that \eqref{eq:u-J} holds. 
Also, since $\phi$ is an analytic function, we have
\begin{equation}
\label{eq:u-I}
u(\phi(I))\subset N\cap (N^-Mh_1Hg\inv).
\end{equation}
Therefore there exists an analytic map $\zeta:I\to M$ such that, 
\[
u(\phi(t))\in N^-\zeta(t)h_1Hg\inv,\quad \forall t\in I, 
\]
and hence
\begin{equation}
\label{eq:limit-phi}
\lim_{\substack{\alpha(a)\to\infty\\a\in A^+}} az(t)u(\phi(t))g\pi(H) 
= z(t)\zeta(t)h_1\pi(H),\quad \forall t\in I.
\end{equation}

Let $\lambda_i|_J$ denote the probability measure as defined in
\eqref{eq:lambda_i} for $J$ in place of $I$.  If there exists a
sequence $j_k\stackrel{k\to\infty}{\longrightarrow}\infty$ such that
$\lambda_{j_k}|_J\stackrel{k\to\infty}{\longrightarrow} \lambda'$,
then by \eqref{eq:limit-phi} we have
\begin{equation}
\label{eq:lambda-prime}
\supp(\lambda')\subset \{z(t)\zeta(t):t\in J\}\cdot h_1\pi(H)\subset\pi(N(H,W)).
\end{equation}
Moreover by Theorem~\ref{thm:W-invariant}, $\lambda'$ is
$W$-invariant. According to the discussion as in
\S\ref{sec:linearization}, we deduce that $\lambda'(\pi(S(H,W))=0$,
because otherwise \eqref{eq:u-J} would hold for a strictly smaller
dimensional subgroup in place of $H$ and some
$g\in\pi\inv(x)$. Therefore, by Proposition~\ref{prop:lambda} and
\eqref{eq:lambda-prime}, for any $f\in\Cc(G/\Gamma)$,
\[
\int_{G/\Gamma} f\,{\dd}\lambda'=\frac{1}{\abs{J}}\int_{t\in J} 
\left(\int_{y\in\pi(H)}
  f(z(t)\zeta(t)h_1y)\,\dd\mu_H(y)\right){\dd}t.
\]
In particular, the right hand side is independent of the choice of the
subsequence $\{j_k\}_{k\in\N}$. Therefore due to
Corollary~\ref{cor:main-measure}, we conclude that for any
$f\in\Cc(G/\Gamma)$,
\begin{equation}
\begin{array}{l}
\lim_{\substack{\alpha(a)\to\infty\\a\in A^+}} \int_{t\in J}
  f(az(t)u(\phi(t))\pi(g))\,{\dd}t  \\
\label{eq:J-mu}
\qquad \qquad = \int_{t\in J}\left(\int_{y\in\pi(H)}
  f(z(t)\zeta(t)h_1y)\,{\dd}\mu_H(y)\right){\dd}t.
\end{array}
\end{equation}

Let $\lambda_t=z_0(t)\zeta(t)\mu_H=z(t)\zeta(t)$ for all $t\in I$,
where $z_0(t)\in M$ such that $z(t)\in z_0(t)A$. We apply
Lemma~\ref{lem:ZA-translate} for $\theta_i(t)=z(t)u(\phi(t))\pi(g)$,
$w_i(t)=z(t)\inv$, and $z_i(t)=e$ for all $t\in I$ and
$i\in\N$. Then \eqref{eq:lambda-muH} follows from
\eqref{eq:J-mu}.
\end{proof}

\begin{proof}[{\bf{Proof of Theorem~\ref{thm:main-unip}}}]
  Let a sequence $g_i\to g$ in $G$ be such that $x_i=\pi(g_i)$. For
  $J$ as in the first paragraph of the proof of Theorem~\ref{thm:psi},
  we consider the discussion of \S\ref{sec:linearization}. By
  Proposition~\ref{prop:lambda} exists a reductive subgroup $H\in\sH$,
  $h_1\in G$, and $\gamma\in\Gamma$ such that $AW\subset h_1Hh_1\inv$
  and by \eqref{eq:n-}
\[
u(\phi(J))\subset N\cap (P^-h_1H(g\gamma)\inv).
\]
Therefore by \eqref{eq:SOn}, \eqref{eq:stereo},
and the analyticity of $\phi$,
\[
S(\phi(I))= p(u(\phi(I)))\subset p(h_1 H\gamma\inv g\inv) = p(H_1h_2),
\]
where $H_1$ is the noncompact simple factor of $h_1Hh_1\inv$
containing $A$ and $h_2=h_1\gamma\inv g\inv\in G$.  In fact, we can
express $H_1=k_1\SO(m,1)k_1\inv$ for some $k_1\in M$ and $2\leq m\leq
n$.  Therefore,
\begin{equation}
  \label{eq:3}
  S(\phi(I))\subset p(\SO(m,1)h_3),  
\end{equation}
where $h_3=k_1\inv h_2$.  Now $p(\SO(m,1))\cong \Sn^{m-1}$, and under
the map $p:G\to \Sn^{n-1}$, the right action of $h_3$ on $G$
corresponds to a conformal transformation on $\Sn^{n-1}$.  Therefore
by \eqref{eq:3}, if $H\neq G$ then $m<n$ and $S(\phi(I))$ is contained
in an $m$-dimensional affine subspace of $\R^{n-1}$ intersecting
$\Sn^{n-1}$. Since $\phi(I)$ is not contained in an affine hyperplane
or a sphere in $\R^{n-1}$ and $S$ is the inverse of stereographic projection,
$S(\phi(I))$ is not contained in a proper subsphere of
$\Sn^{n-1}$. Therefore we conclude that $H=G$.

For each $i\in \N$, we define $\lambda_i|_J$ as in \eqref{eq:lambda_i}
for $J$ in place of $I$. Now if $j_k\to\infty$ is any sequence in $\N$
such that $\lambda_{j_k}|_J\stackrel{k\to\infty}{\longrightarrow}
\lambda_J$ in the space of probability measures, then by the
discussion as in \S\ref{sec:linearization}, by our choice of $H$ as
in Theorem~\ref{thm:Ratner}, and since $H=G$, we have
$\lambda_J(\pi(S(G,W))=0$.  Hence almost all $W$-ergodic components of
$\lambda_J$ are $G$-invariant. Thus $\lambda_J=\mu_G$. Therefore by
Corollary~\ref{cor:main-measure} we conclude that $\lambda_i|_J\to \mu_G$
as $i\to\infty$. Now the conclusion of theorem follows from
\eqref{eq:12} of Lemma~\ref{lem:ZA-translate}.
\end{proof}

\begin{proof}[{\bf {Proof of Theorem~\ref{thm:main}}}] 
  Suppose that \eqref{eq:mainthm} fails to hold for a sequence of
  positive reals $R_i\to\infty$.  Then there exists a sequence
  $\{a_i\}\subset A^+$ such that $\alpha(a_i)\geq R_i$ and a sequence
  $\{x_i\}\subset \cK$ such that
\begin{equation} 
\label{eq:main-fails}
\Abs{\frac{1}{\abs{I}}\int_I f(a_i\theta(t)x_i))\,{\dd}t - \int_{\gmg}
    f\,{\dd}\lambda_G}>\epsilon,\qquad\forall i\in \N.
\end{equation}
Since $\cK$ is compact, by passing to a subsequence, we may assume
that $x_i\to x$ in $\gmg$.

By Bruhat decomposition, $G=P^-N\cup P^-k_0$, where $k_0\in K$ such
that $k_0ak_0\inv=a\inv$ for all $a\in A$.  Also the map $P^-\times
N\to G\setminus\{P^-k_0\}$ given by $P^-\times N\ni (b,u)\mapsto bu$
is an invertible analytic map.  Since $p(\theta(I))$ is not a
singleton set, and $\theta$ is analytic, the set $\{t\in
I:\theta(t)\in P^-k_0\}$ is finite.  As noted earlier, it is enough to
prove the result for all closed subintervals of $I$ with nonempty
interiors and not intersecting this finite set.  Hence without loss of
generality we may assume that $\theta(t)\not\in P^-k_0$ for all $t\in
I$.  Thus we obtain analytic maps $\phi:I\to \R^{n-1}$ and $\zeta:I\to
P^-$ such that
\begin{equation}
  \label{eq:14}
  \theta(t)=\zeta(t)u(\phi(t)),\quad \forall t\in I.
\end{equation}
Then by \eqref{eq:SOn} and \eqref{eq:stereo},
  $S(\phi(I))=p(\theta(I))$.  
  By our assumption, $p(\theta(I))$ is not contained in any hyperplane
  of $\R^{n}$ intersecting $\Sn^{n-1}$. Therefore, since $S$ is the
  inverse of stereographic projection, $\phi(I)$ is not contained in
  any hyperplane or a sphere of $\R^{n-1}$.  Therefore by
  Theorem~\ref{thm:main-unip}, for any subinterval $J\subset I$ with
  nonempty interior,
\begin{equation}
  \label{eq:9a}
  \lim_{i\to\infty} \frac{1}{\abs{J}}\int_J f(a_iu(\phi(t))x_i)\,{\dd}t
  = \int_\gmg f\,{\dd}\mu_G, 
\end{equation}
Now by Lemma~\ref{lem:ZA-translate},
\begin{equation}
  \label{eq:9b}
  \lim_{i\to\infty} \frac{1}{\abs{I}}\int_I
  f(a_i\zeta(t)u(\phi(t))x_i)\,{\dd}t = \int_\gmg f\,{\dd}\mu_G.
\end{equation}
Now \eqref{eq:14} and \eqref{eq:9b} contradict \eqref{eq:main-fails}.
\end{proof}

\begin{proof}[\bf{Proof of Theorem~\ref{thm:Hn}}]
  Let $G=\SO(n,1)$, $K=\SO(n)$, and let $P^-$ be a maximal parabolic
  subgroup of $G$ such that $P^-\cap K=\SO(n-1)$.  Let $A$ the maximal
  $\R$-diagonalizable subgroup of $G$ centralizing $P^-\cap K$. Then
  $A\subset P^-$. Now $G$ admits a transitive right action on
  $\ST(\HH^n)$ via isometries.  We fix $\tilde x_0\in\HH^n$ such that
  $K=\Stab_G(\tilde x_0)$, and we fix $v_0\in S_{\tilde x_0}(\HH^n)$
  such that
\[
\cenKA:=\Stab_K(v_0)=\cen(A)\cap K=\SO(n-1).
\]
Thus
$\ST(\HH^n)\cong \cenKA\backslash G$ and $S_{\tilde x_0}(\HH^n)\cong
\cenKA\backslash K$.  The under this isomorphism, the geodesic flow
$\{\tilde g_t\}$ on $\ST(\HH^n)$ corresponds to the action of
$\{a_t\}=A$ on $\cenKA\backslash G$ by left multiplications, where
$\alpha(a_t)=e^{{\tau}t}$ for all $t\in\R$ and some $\tau>0$.

There exists a discrete subgroup $\Gamma$ of $G$ such that
$\pi:\HH^n\to M$ factors through $\HH^n/\Gamma$ and $M\cong
\HH^n/\Gamma$ as isometric Riemannian manifolds. Hence
$\ST(M)\cong\cenKA\backslash G/\Gamma$, and the geodesic flow $\{g_t\}$ on
$\ST(M)$ corresponds to the left action of $\{a_t\}$ on
$\cenKA\backslash G/\Gamma$.

There exists an analytic map $\theta:I\to G$ such that
\[
\psi(t)=\Der\pi(v_0\theta(t)), \quad \forall t\in I.
\]
As in the proof of Theorem~\ref{thm:main}, let $\phi:I\to \R^{n-1}$ be
the map such that \eqref{eq:14} holds; that is, $\theta(t)\subset
P^-u(\phi(t))$ for all $t\in I$. Then by Theorem~\ref{thm:psi} for
$x=e\Gamma\in G/\Gamma$, there exist $H\in\sH$, $h_1\in G$, and
$\gamma\in\Gamma$ such that $Ah_1\subset h_1H$ and by \eqref{eq:UHg},
$u(\phi(I))\subset P^-h_1H\gamma\inv$. Therefore
\begin{equation}
\theta(t)\subset P^-u(\phi(t))\subset
P^-h_1\gamma\inv=K_0N^-h_1H\gamma\inv, \quad \forall t\in I.
\end{equation}
Therefore, since $\pi(v_0 K_0 g\Gamma)=\pi(v_0g)$, we have
\[
\pi(v_0\theta(t))\subset\pi(v_0 N^-h_1H), \quad A\subset
h_1Hh_1\inv \quad \textrm{and} \quad 
N\cap h_1Hh_1\inv\neq\{e\}. 
\]
Therefore there exists $k_1\in \cenKA$ such that 
\[
\cenKA k_1h_1Hh_1\inv k_1\inv=\cenKA \SO(m,1), \quad
\textrm{where $2\leq m\leq n$}.
\]
Now $v_0\SO(m,1)\cong \ST(\HH^m)$, where $\HH^m$ is isometrically
embedded in $\HH^n$. Since $H\Gamma/\Gamma$ is a closed subset of
$\gmg$,
\[
\pi(v_0h_1H)=\pi(v_0\SO(m,1)k_1h_1)=\pi(\ST(\HH^m)h_2),
\]
is a closed subset of $M$, where $h_2=k_1h_1\in G$.  Therefore $\cenKA
h_1H/(H\cap\Gamma)$ corresponds the embedding of $\Der\Phi(\ST(M_1))$
in $\ST(M)$ which is the derivative of a totally geodesic immersion
$\Phi$ of a hyperbolic manifold $M_1$ in $M$
(see~\cite[\S2]{Shah:tot-geod} for the details). It may also be noted
that the projection of $h_1Hh_1\inv$-invariant probability measure,
say $\mu_1$, on $h_1H/(H\cap\Gamma)$ onto $\cenKA\backslash\gmg\cong
M$, say $\bar\mu_1$, is same as the projection under $\Der\Phi$ of the
normalized measure on $\ST(M_1)$ associated to the Riemannian volume
form on $M_1$.

By \eqref{eq:lambda-muH} of Theorem~\ref{thm:psi}, for any
subinterval $J$ of $I$ with nonempty interior and any
$f\in\Cc(\cenKA \backslash \gmg)$, we have
\begin{equation} 
\lim_{t\to\infty} \frac{1}{\abs{J}}\int_J f(\cenKA a_tu(\phi(s))\Gamma)\,{\dd}s \\
= \int_{\cenKA h_1H\Gamma/\Gamma} f(y)\,{\dd}\bar{\mu}_1(y). 
\end{equation}
Recall that $\theta(s)\in P^-u(\phi(s))$ for all $s\in I$, and $\bar\mu_1$
is $\cen(A)$-invariant with respect to the left action. Therefore by
Lemma~\ref{lem:ZA-translate},
\begin{equation} 
\label{eq:Mgmg}
\lim_{t\to\infty} \frac{1}{\abs{I}}\int_I f(\cenKA a_t\theta(s)\Gamma)\,{\dd}s \\
= \int_{\cenKA h_1H\Gamma/\Gamma} f(y)\,{\dd}\bar{\mu}_1(y). 
\end{equation}

Now in view of the relation between the closed $h_1Hh_1\inv$-orbits
with totally geodesic immersions of finite volume hyperbolic manifolds
as described above, \eqref{eq:Mgmg} implies \eqref{eq:M1-measure}.
\end{proof}

\begin{proof}[\bf {Proof of Corollary~\ref{cor:S}}] 
  Let $\tilde{x}\in \HH^n$ such that $x=\pi(\tilde x)$. We can
  identify $\ST_x(M)$, the unit tangent sphere at $x$, with $\ST_{\tilde
    x}(\HH^n)$, which in turn identifies with the ideal boundary
  sphere $\del\HH^n$ via the visual map.  Since all these
  identifications are conformal, we conclude that
  $\Vis(\tilde\theta(I))$ is not contained in any proper subsphere of
  $\del\HH^n$.  Therefore in terms of the notation in
  Remark~\ref{rem:M1}, $\Sn^{k-1}=\del\HH^n$, and we conclude that
  $M_1=M$ and $\Phi$ is the identity map.  Now the conclusion follows
  from Theorem~\ref{thm:main}.
\end{proof}

\begin{proof}[\bf {Proof of  Theorem~\ref{thm:S-2}}]
  The proof is similar to the proof of Corollary~\ref{cor:S}.
\end{proof}

\begin{proof}[\bf {Proof of Theorem~\ref{thm:Hn-basic}}]
  We identify $\Sn^{n-1}$ with a hyperbolic sphere of radius $1$
  centered at $0$ in $\HH^n$ (in the unit Ball $B^n$-model), say $S$,
  and treat $\bar\psi$ as a map from $I$ to $S$. For any $s\in I$, let
  $v_s\in \ST_{\bar\psi(s)}(\HH^n)$ be the unit vector normal to $S$
  which is also a tangent to the directed geodesic from $0$ to
  $\bar\psi(s)$. We define an analytic curve $\psi:I\to \ST(M)$ by
\[
\psi(s)=(\pi(\bar\psi(s)),\Der\pi(v_s)), \quad \forall s\in I.
\]
Therefore the condition of Theorem~\ref{thm:S-2} is satisfied, because
\[
\Vis(\tilde\psi(s))=\Vis((\bar\psi(s),v_s))=\bar\psi(s), \forall s\in I,
\]
and hence $\Vis(\tilde\psi(s))$ is not contained in a proper subsphere
of $\del\HH^n$. 

For any $\alpha>0$, $\pi(\alpha\bar\psi(s)) =
g_{t(\alpha)}\pi(\psi(s))$ for some $t(\alpha)>0$ such that
$t(\alpha)\to\infty$ as $\alpha\to 1^-$.  Therefore \eqref{eq:Hn}
follows from Theorem~\ref{thm:S-2}.
\end{proof}

\subsection{A stronger version}

\begin{theo} 
  \label{thm:G} Let $G=\SO(n,1)$ and $\Gamma$ a lattice in $G$.  Let
  $\theta:I=[a,b]\to G$, where $a<b$, be an analytic map such that for
  any minimal parabolic subgroup $P^-$ of $G$, the image of $\psi(I)$
  in $P^-\backslash G\cong \Sn^{n-1}$ is not contained in any proper
  subsphere of $\Sn^{n-1}$. Then given any $f\in \Cc(\gmg)$, a
  compact set $\cK\subset G$ and an $\epsilon>0$, there exists a
  compact set $\cC\subset G$ such that
  \[
  \Abs{\frac{1}{\abs{I}}\int_I f(g\theta(t)x)\,{\dd}t - \int_{\gmg}
    f\,{\dd}\lambda_G}<\epsilon,\quad \forall x\in\cK \textrm{ and
  }\forall g\in G\setminus\cC.
  \]
\end{theo}

A proof of the above generalization of Theorem~\ref{thm:main} can be
given by similar arguments. The analogue of \S\ref{subsec:conseq} is a
little more delicate in this case. We do not include the proof here in
order to have simpler proofs for all other results.

\section{Scope for generalization and applications} 
 
The results of this article lead to obvious similar questions about
expanding translates of $(C,\alpha)$-good curves on horospherical
subgroups of general semisimple Lie groups. Especially the
affirmative answer to the following question has interesting
applications to problems in Diophantine approximation
\cite{Klein+Mar:Annals98,Kleinbock+Weiss:Dirichlet}:

\begin{quest} \label{quest:sln}
 Let $G=\SL(n+1,\R)$, $\Gamma=\SL(n+1,\Z)$, and $\mu_G$
  denote the $G$-invariant probability measure on $\gmg$.  Let
  \begin{equation*}
    u(v)=
    \left(\begin{smallmatrix} 
      1  & v \\
      0  & I_{n} 
    \end{smallmatrix}\right),\ \forall v\in \R^n,  \textrm{ and }
    a(t)=\diag(e^{nt},e^{-t},\dots,e^{-t}),\  \forall t\in\R.  
  \end{equation*}
  Let $\phi:[0,1]\to \R^n$ be an analytic (or a $(C,\alpha)$-good)
  curve such that its image is not contained in any proper affine
  hyperplane in $\R^n$.  Then is it true that for any $x\in\gmg$ and
  any $f\in\Cc(\gmg)$,
  \begin{equation}
    \label{eq:99}
    \lim_{t\to\infty} \int_0^1f(a(t)u(\phi(s))x)\,{\dd}s=\int_\gmg f\,{\dd}\mu_G\,?
  \end{equation}
\end{quest}

The main result of \cite{Klein+Mar:Annals98} provides a very good
estimate on the rate of nondivergence of this translated measure.  The
method of this article is applicable to show that, after a suitable
modification of the curve by elements form the centralizer of
$\{a(t)\}$, the limiting measure is invariant under a unipotent
one-parameter subgroup of the form $\{u(sw_0)\}$ for some
$w_0\in\R^n\setminus\{0\}$.  Also the method to study behaviour of
expanded trajectories near the singular sets is applicable
here. Obtaining an analogue of Lemma~\ref{lem:SL2} in order to derive
algebraic consequences of Proposition~\ref{prop:main3} is the main
difficulty in this problem.

Since the initial submission of this article, the
author~\cite{Shah:SLn} has answered Question~\ref{quest:sln} in
affirmation for analytic curves.

In another direction, it is still an open question to prove the exact
analogue of Theorem~\ref{thm:S-2} for the actions of $\SO(n,1)$ on
homogeneous spaces of larger Lie group $G$ containing $\SO(n,1)$; see
\cite{Goro:AIM}.

\end{document}